\documentclass[]{informs3}
\OneAndAHalfSpacedXI

\usepackage{amsmath,amsfonts,amssymb}
\usepackage{url}
\usepackage[capitalize,noabbrev]{cleveref}
\usepackage{xcolor}
\usepackage{nicefrac}
\usepackage{booktabs}
\usepackage{anyfontsize}
\usepackage{siunitx}
\usepackage[title]{appendix}
\usepackage{endnotes}
\def\notesname{Endnotes}%
\def\makeenmark{\hbox to1.275em{\theenmark.\enskip\hss}}
\def\enoteformat{\rightskip0pt\leftskip0pt\parindent=1.275em
\leavevmode\llap{\makeenmark}}
	
\newcommand{\blue}[1]{{\color{black}#1}}

\colorlet{darkorange}{orange!80!black}
\newcommand{\refone}[1]{{\color{black}#1}}
\newcommand{\reftwo}[1]{{\color{black}#1}}
\newcommand{\rev}[1]{{\color{black}#1}}

\newcommand{\R}{\mathbb{R}_+}
\newcommand{\B}{\{0,1\}}
\newcommand{\BL}{[0,1]}

\newcommand{\aop}{AOP\xspace}
\newcommand{\probname}{AOPC\xspace}
\newcommand{\pmin}{p_{\text{min}}}
\newcommand{\pmax}{p_{\text{max}}}
\newcommand{\vx}{v(\gb{x})}
\newcommand{\unp}{\un{p}}
\newcommand{\ovp}{\ov{p}}
\newcommand{\ovzp}{\ov{z}(\unp,\ovp)}

\newcommand{\ovzpk}{\ov{z}(p^k,p^{k+1})}
\newcommand{\unzp}{\un{z}(\unp,\ovp)}
\newcommand{\unzpk}{\un{z}(p^k,p^{k+1})}
\newcommand{\gexp}{G^{\text{exp}}_\rho}

\newcommand{\gb}[1]{\boldsymbol{#1}}
\newcommand{\st}{\text{s.t.}\;\;}
\newcommand{\ov}[1]{\overline{#1}}
\newcommand{\un}[1]{\underline{#1}}

\newcommand{\tabsp}{\addlinespace[0.5em]}
\newcommand{\cpuavg}{cpu$_{\text{avg}}$}
\newcommand{\cpumax}{cpu$_{\text{max}}$}
\newcommand{\our}{LLRS}
\newcommand{\intavg}{int$_{\text{avg}}$}
\newcommand{\intmax}{int$_{\text{max}}$}
\newcommand{\comb}{(n,\varPhi,\gamma)}
\newcommand{\final}{MILP+}
\newcommand{\rhof}{{$\rho^f~$}}
\newcommand{\rhol}{{$\rho^\ell~$}}

\let\emptyset\varnothing

\usepackage{natbib}
\bibpunct[, ]{(}{)}{,}{a}{}{,}%
\def\bibfont{\small}%

\TheoremsNumberedThrough
\ECRepeatTheorems
\EquationsNumberedThrough 

\newcommand{\eof}{\hfill$\square$}

\usepackage[linesnumbered,ruled,vlined]{algorithm2e}
\DontPrintSemicolon
\SetKw{KwAnd}{and}
\SetKw{KwOr}{or}
\SetKw{Break}{break}

\begin{document}
	
\RUNAUTHOR{Leitner et al.}
\RUNTITLE{Exact assortment optimization with product costs}
\TITLE{An Exact Method for (Constrained) Assortment Optimization Problems with Product Costs}
	
\ARTICLEAUTHORS{
	\AUTHOR{Markus Leitner}
	\AFF{Department of Operations Analytics, Vrije Universiteit Amsterdam, The Netherlands\\ \EMAIL{m.leitner@vu.nl}}
	\AUTHOR{Andrea Lodi}
	\AFF{Jacobs Technion-Cornell Institute, Cornell Tech and Technion - IIT, USA\\
	\EMAIL{andrea.lodi@cornell.edu}}
	\AUTHOR{Roberto Roberti}
	\AFF{Department of Information Engineering, University of Padova, Italy\\ \EMAIL{roberto.roberti@unipd.it}}		
	\AUTHOR{Claudio Sole}
	\AFF{Canada Excellence Research Chair in Data-Science for Real-time Decision-Making, Polytechnique Montr\'eal, Canada\\ 
	\EMAIL{claudio.sole@polymtl.ca}}
}

\ABSTRACT{
We study the problem of optimizing assortment decisions in the presence of product-specific costs when customers choose according to a multinomial logit model. This problem is NP-hard and approximate solutions methods have been proposed in the literature to obtain both \refone{lower} and \refone{upper} bounds in a tractable manner. We propose the first exact solution method for this problem and show that provably optimal assortments of instances with up to one thousand products can be found, on average, in about two tenths of a second. In particular, we propose a bounding procedure based on the approximation method of \cite{Feldman2015} to provide tight \refone{lower} and \refone{upper} bounds at a fraction of their computing times. We show how these bounds can be used to effectively identify an optimal assortment. We also describe how to adapt our approach to handle {cardinality or space/resource capacity constraints on the assortment as well as assortment optimization under a mixed-multinomial logit model. In both cases, our solution method provides significant computational boosts compared to exact methods from the literature.}
}
	
	
\maketitle

\section{Introduction\label{sec:intro}}
Composing a set of products that maximizes \refone{the expected profit} is a major planning problem in retail operations, revenue management, and online advertising, see, e.g., \citet{Gallego2019,Koek2008}. The need to consider customer choice behavior in this decision has led to the consideration of a large variety of corresponding models, particularly including parametric discrete-choice models based on random utility theory \citep{Strauss2018}. This led to a vast body of literature on \textit{assortment optimization problems} (\aop) under variants and extensions of the \textit{multinomial logit model} (MNL) \refone{(see, e.g., \cite{Gallego2019})}. Despite some drawbacks, such as the assumption of independence of irrelevant alternatives that are overcome by some of its extensions, the MNL remains one of the most frequently studied discrete-choice models in assortment optimization. \refone{The MNL is so popular for two main reasons: ($i$) it can be efficiently estimated, and ($ii$) {choice probabilities under this model take closed form and are easily interpretable. We refer again the interested reader to the works of \cite{Koek2008}, \cite{Strauss2018} and references therein for a review of the applications of this model in fields such as revenue management, marketing and travel industry among others.}}

We focus on variants of the \textit{assortment optimization problem with product costs} (\probname) under the MNL as introduced by \citet{Kunnumkal2009} (see, also, \citet{Kunnumkal2019}) who also show that the \probname is NP-hard. The \probname is particularly important as it models a variety of practical applications where a certain cost (e.g., stocking 
costs) is incurred by the firm for offering a given product. Furthermore, the \probname appears as a subproblem when solving other problems encountered in assortment optimization and revenue management where business constraints, such as space \citep[see, e.g.,][]{Feldman2019} or resource capacity constraints \citep[see, e.g.,][]{kunnumkal2010newDecomp}, must be satisfied, or when optimizing assortment decisions over a mixture of customer types \citep[see, e.g.,][]{Feldman2015}. In particular, (approximate) solutions schemes based on the Lagrangian relaxation that rely on the solution of AOPs with product costs have been proposed to solve these problems \citep[see, e.g.,][]{Feldman2015}.

Building upon the parametric \refone{upper} bound presented by \cite{Feldman2015} and a mixed-integer linear programming (MILP) \rev{model} presented by \cite{Kunnumkal2019}, we propose an exact method that can find an optimal assortment for \probname instances featuring up to 1000 products in about two tenths of a second of computing time on average. To the best of our knowledge, this is the first exact algorithm that is not based on simply casting a MILP model or a mixed-integer conic quadratic (MICQ) model into a general-purpose black-box solver. We show that our method significantly outperforms these models in terms of size of the instances that can be solved and in terms of computing time by orders of magnitude.

The remainder of this article is organized as follows. \cref{sec:litrev} summarizes the literature related to the \probname. \cref{sec:prob} formally introduces the \probname and reviews two exact formulations from the literature that can be solved with general-purpose black-box solvers. \cref{sec:framework} introduces the exact solution method we propose. \cref{sec:constrained} shows how the method can be extended to the constrained \probname. Computational results are reported in \cref{sec:results}.
{\cref{sec:mmnl} discusses adaptations required to apply the proposed solution algorithm to assortment optimization under mixtures of multinomial logits and reports computational results obtained for this problem variant.}
Some conclusions are drawn in \cref{sec:conclusions}.

\section{Literature review\label{sec:litrev}}
In the following, we provide a brief literature review mainly focusing on assortment optimization problems under variants of the MNL that are directly relevant for our developments. 

\citet{Talluri2004} show that the \aop under the MNL can be solved in polynomial time by considering revenue-ordered product subsets. While the capacity constrained \aop (under the MNL) can still be solved in polynomial time \citep{Rusmevichientong2010}, most other practically interesting problem variants are NP-hard. This includes, in particular, the \probname \citep{Kunnumkal2009} and the \aop under the more general mixed-multinomial logit model (MMNL), which considers multiple customer classes \citep{Bront2009, Rusmevichientong2014}.
Consequently, a large body of literature is devoted to the study of (polynomial time) approximation algorithms for AOPs, see, e.g., \citet{Feldman2015a,Feldman2019a,Gallego2014,Liu2020} for recent contributions.

Despite being a natural and important generalization, the \probname has been addressed by few articles only. \citet{Kunnumkal2009} and \citet{Kunnumkal2019} reformulate the \probname as a parametric optimization problem and show how to derive \refone{lower} and \refone{upper} bounds
on the optimal expected profit in polynomial time. The former work proposes a 2-approximation algorithm with complexity $O(n^3)$, where $n$ is number of products considered, and a polynomial time approximation scheme for computing feasible assortments to the \probname with theoretical guarantees on their expected revenue.  Complementary to this work, \cite{Kunnumkal2019} propose a procedure to obtain 
upper bounds
on the optimal revenue that relies on the solution of $O(n^3)$ continuous knapsack problems. 
Their approach favorably compares against a MILP formulation for the \probname provided by the authors in the same work (which we will recall in Section \ref{sec:prob}), with relatively tight bounds obtained in short computing times.

The \probname has also been considered in other AOP variants. \citet{Feldman2015} suggest a solution method based on 
Lagrangian relaxation to approximately solve the \aop under the MMNL. They relax the constraints ensuring that the same products are offered to all customer types and therefore obtain the \probname as subproblem. In order to obtain \refone{upper} bounds to the \probname in a tractable manner, they propose a grid-based approach where denser grids 
provide tighter bounds, at the cost of higher computing times. \citet{Feldman2019} relate the approximability of the space-constrained AOP to the solution of AOPs with fixed costs. \citet{Honhon2012} propose polynomial time algorithms for 
different variants of an AOP with fixed costs under a ranking-based customer choice model. In the context of network revenue management, the column-generation procedure proposed by \cite{kunnumkal2008refined} relies on the solution of a subproblem 
that has 
the form of an \probname, and \cite{kunnumkal2010newDecomp} propose a Lagrangian-decomposition approach that relies on the solution of a series of AOPs with fixed costs, where the Lagrangian multipliers associated with the relaxed constraints play the role of product fixed costs.

As mentioned above, the literature on exact methods for NP-hard AOPs under variants of the MNL is relatively limited. 
In this regard, \citet{Bront2009} and \cite{Mendez-Diaz2014} show that an \aop under the (M)MNL can be reformulated as a MILP with a polynomial number of variables and constraints. Although solving this compact MILP (by a black-box solution framework) does not scale well in practice, such approach is typically used to assess the performance of alternative methods. \citet{Sen2018} suggest to reformulate the (constrained) AOP under the MMNL as a conic (quadratic) mixed-integer program 
that can be solved by black-box solution frameworks too.
The results of their computational study indicate 
clear advantages of their approach over previously considered MILP-based methods. Their formulation will be adapted to the \probname in \cref{sec:prob} and considered in our computational study.  
%
Recently, \citet{Alfandari2021} propose an exact algorithm based on fractional programming for the \aop under the nested logit model. 

The paper that is closest to ours is the one by \cite{Feldman2015}, which we will describe in Section \ref{subsec:ft}. In particular, one of the contributions of our work is to embed their grid-based approximation method into a bounding procedure, where coarser grids are used to alleviate the computational burden of denser ones, allowing to obtain both \refone{lower} and \refone{upper} bounds of the same quality at a small fraction of their computing times. Furthermore, we show how to leverage such bounds to solve the \probname to optimality.

\section{Problem description and formulations from the literature\label{sec:prob}}
The \probname can be formally described as follows. A set of $n$ products $P = \{1,2, \dots, n\}$ is given. Each product $j \in P$ is characterized by a revenue $r_j > 0$, a product {cost} $c_j\ge 0$, and a preference weight (or, simply, preference) $v_j>0$. The preference of not making any purchase is denoted as $v_0 \geq 0$. Let $x_j\in \{0,1\}$ be a binary decision variable indicating if product $j\in P$ is included in the assortment ($x_j = 1$) or not ($x_j = 0$). For each $\gb{x} \in \B^n$, let $\vx$ be defined as $\vx = \sum_{j \in P} v_j x_j$. According to the MNL, the probability that a customer purchases product $j \in P$ is $v_j x_j / (v_0 + \vx)$, and the no-purchase probability is $v_0 / (v_0 + \vx)$. The objective of the \probname is to find a set of products $P^\ast \subseteq P$ such that the corresponding profit $z(P^\ast) = \sum_{j \in P^\ast} r_j v_j / ( v_0 + \sum_{j \in P^\ast} v_j) - \sum_{j\in P^\ast} c_j$ is maximum. The \probname can be formulated as the 
mixed integer non-linear programming problem
\begin{align}
z^\ast = \max_{\gb{x} \in \B^n} \bigg\{ \frac{\sum_{j\in P} r_j v_j x_j}{v_0 + \vx} - \sum_{j\in P} c_j x_j \bigg\}.
\label{prob}
\end{align}

Two exact reformulations of problem \eqref{prob} that can be cast into general-purpose black-box solvers have been proposed recently. The first formulation is a MILP model provided by \citet{Kunnumkal2019}, see also \citet{Gallego2019}. For each $j \in P$, in addition to decision variable $x_j \in \B$, let $u_j \geq \R$ be a non-negative continuous variable representing the purchasing probability of product $j \in P$, i.e., $u_j = v_j x_j / (v_0 + \vx)$. Similarly, let $u_0\geq \R$ be a continuous variable representing the no-purchase probability, i.e., $u_0 = v_0 / (v_0 + \vx)$. The MILP provided by \citet{Kunnumkal2019} is
\begin{subequations}\label{prob:milp}
	\begin{align}
	z^\ast = & \max \;\; \sum_{j \in P} \left( r_j u_j - c_j x_j \right) \label{prob:milp:obj} \\
	\st & v_0 u_j \leq v_j u_0 & \forall j\in P \label{prob:milp:prob:ratios} \\
	& u_j \leq \frac{v_j}{v_0+v_j} x_j & \forall j\in P \label{prob:milp:force} \\
	& u_0 + \sum_{j \in P} u_j = 1 \label{prob:milp:conv} \\
	& x_j \in \B & \forall j \in P \label{prob:milp:int} \\
	& u_j \in \R & \forall j \in P \cup \{0\} \label{prob:milp:cont}
	\end{align}
\end{subequations}
The objective function \eqref{prob:milp:obj} maximizes the profit of the selected products. Constraints \eqref{prob:milp:prob:ratios} ensure that the ratio between the purchasing probabilities of product $j \in P$ and the no-purchase probability is consistent with the corresponding preferences $v_j$ and $v_0$. Constraints \eqref{prob:milp:force} force the purchasing probabilities of all products not included in the assortment to be zero. 
Constraint \eqref{prob:milp:conv} makes sure that the sum of purchasing probabilities is equal to one. Constraints \eqref{prob:milp:int}-\eqref{prob:milp:cont} define the \rev{domains} of the variables.

The second formulation is based on a MICQ proposed by \citet{Sen2018} for the AOP under the MMNL. In particular, we adjust the corresponding objective function 
 to the \probname by adding a costs-related term, which penalizes the introduction of products in the assortment. Specifically, let $\ov{r}$ be the maximum revenue over all products, i.e., $\ov{r} = \max \{ r_j \, | \, j \in P \}$.  \rev{Let} $\varphi_j \in \R$ be a non-negative continuous variable equal to $1 / (v_0 + \vx)$ if product $j \in P \cup \{0\}$ is in the assortment and 0 otherwise - notice that $\varphi_0$ is equal to $1 / (v_0 + \vx)$. Moreover, let $w \in \R$ be a non-negative continuous 
variable equal to the sum of the preferences of the selected products including the no-purchase preference. The \probname\ can be formulated as
\begin{subequations} \label{SAK}
	\begin{align}
	z^\ast = \;\; & \ov{r} - \min  \big( \ov{r} v_0 \varphi_0 + \sum_{j \in P} v_j ( \ov{r} - r_j ) \varphi_j + \sum_{j \in P} c_j x_j \big) \label{SAK0} \\
	\st & w = v_0 + \vx & \label{SAK1} \\
	& \varphi_j w \geq x_j^2 & \forall j\in P \label{SAK2} \\
	& \varphi_0 w \geq 1 & \label{SAK3} \\
	& v_0 \varphi_0 + \sum_{j \in P} v_j \varphi_j \geq 1 \label{SAK4} \\
	& x_j \in \B & \forall j \in P \label{SAKa} \\
    & \varphi_j \in \R  & \forall j \in P \cup \{0\} \label{SAKb} \\
    & w \in \R & \label{SAKc}
	\end{align}
\end{subequations}


The objective function \eqref{SAK0} aims at maximizing (or equivalently, minimizing the negative) expected profit. Here, the first three terms denote the revenue of the selected products, while the last term  accounts for the total cost of adding \rev{them}
to the assortment. 
Constraint \eqref{SAK1} sets variable $w$ equal to the sum of the preferences of the selected products plus the no-purchase preference. Constraints \eqref{SAK2} guarantee that $\varphi_j = 1 / w$ if product $j \in P$ is selected and 0 otherwise, since the problem is in minimization form and variables $\varphi_j$ have non-negative objective coefficients. Constraint \eqref{SAK3} sets $\varphi_0$ equal to $1 / w$. Constraint \eqref{SAK4} corresponds to \eqref{prob:milp:conv} and is redundant but helps significantly strengthen the \rev{continuous} 
relaxation of \eqref{SAK} (see \citet{Sen2018}). Constraints \eqref{SAKa}-\eqref{SAKc} define the domain of the variables.

\citet{Sen2018} also show how to strengthen the \rev{continuous}
relaxation of formulation \eqref{SAK} by adding the following McCormick inequalities \citep{mccormick1976}:
\begin{subequations} \label{McC}
	\begin{align}
	& m^1_j x_j \leq \varphi_j \leq M^1_j x_j & \forall j \in P \label{McC1} \\
	& \varphi_0 - M^0_j ( 1 - x_j) \leq \varphi_j \leq \varphi_0 - m^0_j (1 - x_j) & \forall j \in P \label{McC2}
	\end{align}
\end{subequations}
where $m^1_j$ and $M^1_j$ are appropriate lower and upper bounds on the value of variable $\varphi_j$ when product $j \in P$ is selected, $M^0_j$ and $m^0_j$ are an upper and a lower bound on the value of variable $\varphi_0$ when product $j \in P$ is not selected. For the unconstrained \probname, $m^1_j$ and $M^1_j$ are set equal to $1 / (v_0 + \sum_{i \in P} v_i)$ and $1 / (v_0 + v_j)$, respectively, whereas $M^0_j$ and $m^0_j$ are set equal to $1 / (v_0 + \sum_{i \in P \setminus \{j\}} v_i)$ and $1 / v_0$, respectively. For the constrained \probname (i.e., when cardinality or budget constraints are present), tighter values of parameters $m^1_j$, $M^1_j$, $M^0_j$, and $m^0_j$ can be computed by solving some binary (multiple) knapsack problems or the corresponding linear relaxations (see \citet{Sen2018}).

\section{Exact solution method\label{sec:framework}}
In this section, we describe our exact method, which builds upon the parametric \refone{upper} bound proposed by \cite{Feldman2015} (summarized in \cref{subsec:ft}). The exact method is based on three main ideas. The first idea (see \cref{subsec:seqbound}) is to enhance the bounding procedure of \citet{Feldman2015} by limiting the number of parametric bounds to compute and sequentially improve the global \refone{upper} bound with the goal of finding a small range for the sum of the preferences of the products of any \reftwo{optimal \probname assortment}. The second idea (see \cref{subsec:varfix}) is to use the \refone{lower} and \refone{upper} bounds computed by the bounding procedure to rule out the products that cannot be part of any \reftwo{optimal \probname assortment}. The third idea (see \cref{subsec:gap}) is to find an \reftwo{optimal \probname assortment} by solving problem \eqref{prob:milp}, with a general-purpose solver, and add some constraints to limit the search within the range of preferences returned by the bounding procedure.

\subsection{The approximation method of \citet{Feldman2015}} \label{subsec:ft}
\citet{Feldman2015} propose the following approximation method that computes (provably) tight \refone{upper} bounds to the \probname. Their method is based on formulating the \probname as
\begin{equation}
z^\ast = \max_{p \in [\pmin,1]} \{ z(p) \} \label{FTzstar}
\end{equation}
where $p$ is the no-purchase probability and $\pmin$ is the minimum no-purchase probability (achieved when all products are in the assortment). By normalizing the purchase preferences such that $v_0=1$ (which can be done w.l.o.g.), these values are computed as $p=1 / (1 + \vx)$ and $\pmin = 1 / (1 + \sum_{j \in P} v_j)$. Furthermore, $z(p)$ is the optimal value of the parametric problem
\begin{equation}
z(p) = \max_{\gb{x} \in \B^n} \Bigg\{ \sum_{j \in P} \left( {p \; r_j} v_j - c_j \right)x_j \, | \, \frac{1}{1 + \vx} = p \Bigg\}, \label{FTpar}
\end{equation}
which can be equivalently written as (see \citet{Feldman2015})
$$
z(p) = \max_{\gb{x} \in \B^n} \Bigg\{ \sum_{j \in P} \left( {p \; r_j} v_j - c_j \right)x_j \, | \, \vx \leq \frac{1}{p} - 1 \Bigg\}.
$$

Solving \eqref{FTzstar} by computing $z(p)$ is intractable because it requires solving a binary knapsack problem, which is NP-hard, for each value of $p \in [\pmin,1]$. The parametric problem \eqref{FTpar} can, however, be adjusted to compute provably tight \refone{upper} bounds to $z^\ast$ as follows.

For any pair of values $\unp,\ovp \in [\pmin,1]$ such that $\unp \leq \ovp$, let $\ovzp$ be a parametric \refone{upper} bound to $\max_{p \in [\unp,\ovp]} z(p)$ (i.e., $\ovzp \geq \max_{p \in [\unp,\ovp]} z(p)$) defined as
\begin{equation}
\ovzp = \max_{\gb{x} \in \BL^n} \Bigg\{ \sum_{j \in P} \left( {\ovp \; r_j} v_j - c_j \right)x_j \, | \, \vx \leq \frac{1}{\unp} - 1 \Bigg\}. \label{FTparub}
\end{equation}

\citet{Feldman2015} show that 
\begin{equation}
z(G) = \max_{k \in \{1,\ldots,K\}} \big\{ \ovzpk \big\} \label{FTubg}
\end{equation}
is a valid \refone{upper} bound to $z^*$ for any set of grid points $G = \{ p^k \, | \, k \in \{ 1, \ldots, K+1 \} \}$, such that $\pmin = p^1 \leq p^2 \leq \ldots \leq p^K \leq p^{K+1} = 1$, which can be computed by solving $K$ continuous knapsack problems. Each of these knapsack problems can be solved by considering and adding (if feasible) the products $j\in P$ in non-increasing order of their utility-to-space consumption ratios $(\ovp r_jv_j - c_j) / v_j$.

The gap between $z(G)$ and $z^\ast$ is due to two sources of \textit{error}. The first source is that $z^\ast$ is computed over all values of $p$ in the interval $[\pmin,1]$ whereas $z(G)$ is computed over a set of grid points and $p$, in \eqref{FTpar}, is replaced by $\unp$ and $\ovp$ in the objective function and the right-hand side of the constraint of \eqref{FTparub}, respectively. Hence, denser grids provide better approximations of $z^\ast$. The second source is that the decision variables to compute $\ovzp$ are in the range $\BL^n$ whereas the \probname imposes the integrality constraints on each $x_j$, $j \in P$; nevertheless, the linear relaxation of the binary knapsack problem usually provides tight \refone{upper} bounds.

\citet{Feldman2015} discuss properties of effective grids and illustrate the benefits of using an exponential grid. For a fixed parameter $\rho > 0$ (e.g., $\rho = 0.1, 0.01, 0.001, \ldots$), the exponential grid $\gexp$ is defined as $\gexp = \{ (1 + \rho)^{-k+1} \, | \, k \in \{ 1,\ldots,K(\rho)+1 \} \}$, where $K(\rho)$ is such that $(1 + \rho)^{-K(\rho)} \leq \pmin < (1 + \rho)^{-K(\rho)+1}$. The exponential grid of points $\gexp$ clearly covers the entire interval $[\pmin,1]$. \citet{Feldman2015} show that $z(\gexp)$ cannot be improved by more than a factor $1 + \rho$ by any other grid. Since $\pmin < (1 + \rho)^{-K(\rho)+1}$, we have $K(\rho) = O(- \log(\pmin) / \log(1 + \rho) )$. Therefore, for example, if $\pmin = 0.25$ and we would like to have a performance guarantee of $0.1\%$, we can choose $\rho = 0.001$, which implies $K(\rho) = 1387$. In other words, tight \refone{upper} bounds to $z^\ast$ can be achieved by using exponential grids with a relatively small number of points. 

\subsection{Enhanced bounding procedure} \label{subsec:seqbound}
We now discuss how $z(\gexp)$ can be computed also for very dense exponential grids by sequentially computing \refone{upper} bounds over grids with increasing density.

For each pair of values $\unp, \ovp \in [\pmin,1]$ such that $\unp \leq \ovp$, let $\unzp$ be a \refone{lower} bound to $z^\ast$. We compute this \refone{lower} bound by considering the products $j\in P$ in non-increasing order of their utility-to-space consumption ratio $(\ov{p}r_jv_j - c_j) / v_j$ and adding a given product to the assortment $P'\subseteq P$ as long as the sum of the preferences of the products in $P'$ does not exceed $\frac{1}{\un{p}} - 1$. This simple heuristic for the binary knapsack problem, which usually provides tight \refone{lower} bounds, yields an assortment $P'$ with a profit of $\unzp=\sum_{j \in P'} (pr_jv_j - c_j)$, where $p = 1 / (1 + \sum_{j \in P'} v_j)$. 

For a given $\gexp$ and the corresponding parametric bounds $\unzpk$, $\ovzpk$ (for $k \in \{ 1,\ldots,K(\rho) \}$), let
\begin{equation}
lb(\gexp) = \max_{k \in \{1,\ldots,K(\rho)\}} \{ \unzpk \} 
\end{equation}
be the best \refone{lower} bound to $z^\ast$ computed over $\gexp$. Moreover, let $\mathcal{I}(\gexp)$ be the set of intervals, defined by $\gexp$, that contains a value $p$ for which $z(p)$ is {a candidate} optimal solution of \eqref{FTzstar}. Set $\mathcal{I}(\gexp)$ is defined as
\begin{equation}
\mathcal{I}(\gexp) = \{ [p^k,p^{k+1}] \; | \; k \in \{1,\ldots,K(\rho)\} \, : \, \ovzpk \geq lb(\gexp) \}.
\end{equation}
Let $\mathcal{P}(\gexp)$ be the union of all intervals of $\mathcal{I}(\gexp)$, i.e., $\mathcal{P}(\gexp ) = \cup_{[p^k,p^{k+1}] \in \mathcal{I}(\gexp)} [p^k,p^{k+1}]$. 
We will 
refer to the minimum and maximum values of $\mathcal{P}(\gexp)$ as $\pmin(\gexp)$ and $\pmax(\gexp)$, respectively, 
i.e., $\pmin(\gexp) = \min \{ p \, | \, p \in \mathcal{P}(\gexp) \}$ and $\pmax(\gexp) = \max \{ p \, | \, p \in \mathcal{P}(\gexp) \}$.

As $lb(\gexp)$ is a \refone{lower} bound to $z^\ast$ and $\ovzpk$ is an \refone{upper} bound to $z(p)$ with $p \in [p^k,p^{k+1}]$, $k\in \{1, \dots, K(\rho)\}$, we can observe that, for any optimal \probname solution $P^\ast$, there must exist an interval $[\unp,\ovp] \in \mathcal{I}(\gexp)$ such that $p^\ast \in [\unp,\ovp]$, where $p^\ast = 1 / (1 + \sum_{j \in P^\ast} v_j)$, so $p^\ast$ belongs to $\mathcal{P}(\gexp)$.

This observation suggests that we can compute $z(\gexp)$ even for highly dense grids (e.g., defined with $\rho=\num{1E-07}$), by exploiting sets $\mathcal{I}(\gexp)$ computed with much higher values of $\rho$ (e.g., $\rho=\num{1E-2}$ or $\num{1E-3}$). In particular, given $\rho'$ and $\rho$ such that \reftwo{$\rho' > \rho$}, $z(\gexp)$ can be computed as 
\begin{equation}
z(\gexp) = \max_{k \in \{1,\ldots,K(\rho)\} \, : \, [p^k,p^{k+1}] \cap \mathcal{P}(G^{\text{exp}}_{\rho'}) \neq \emptyset} \{ \ovzpk \}. \label{ourpub}
\end{equation}

We propose an \refone{upper} bounding procedure that computes a sequence of \refone{upper} bounds to $z^\ast$ of increasing tightness {using two parameters $\rho^f$ and $\rho^\ell$. 
The procedure first computes $z(\gexp)$ according to \eqref{FTubg} with $\rho = \rho^f$. Subsequently, it computes $z(\gexp)$ according to \eqref{ourpub} for $\rho = 0.1\rho^f, 0.01\rho^f, \ldots, \rho^\ell$ and by using $\mathcal{P}(G^{\text{exp}}_{\rho'})$ with $\rho' = \rho^f, 0.1\rho^f, \ldots$.
Our computational results show that the number of intervals for which the parametric bound \eqref{FTparub} must be computed with this procedure is significantly lower than by computing $z(\gexp)$ as in \eqref{FTubg}, with $\rho = \rho^\ell$.
}


\subsection{Variable fixing to rule out products} \label{subsec:varfix}
The bounding procedure described in \cref{subsec:seqbound} can be \reftwo{sped} up by ruling out some products that cannot be part of any \reftwo{optimal \probname assortment}. For each value of $\rho = \rho^f, 0.1\rho^f, \ldots, \rho^\ell$, the bounding procedure returns \refone{lower} and \refone{upper} bounds ($lb(\gexp)$ and $z(\gexp)$\rev{, respectively}) to $z^\ast$ and the set $\mathcal{P}(\gexp)$, which allow to fix the variable $x_j$ of some products $j\in P$ to 0 according to the following propositions.

\begin{proposition}
For any $\rho$, each product $j \in P$ such that $\pmax(\gexp) r_j v_j - c_j < 0$ cannot be part of any \reftwo{optimal \probname assortment}.
\end{proposition}
\textit{Proof.} By definition, $\pmax(\gexp)$ is greater than or equal to the no-purchase probability of any \reftwo{optimal \probname assortment}. If $\pmax(\gexp) r_j v_j - c_j < 0$ for a given product $j \in P$, then $p r_j v_j - c_j$ is negative for any $p \in \mathcal{P}(\gexp)$. Therefore, such product cannot yield a positive profit and cannot be part of any \reftwo{optimal \probname assortment}.\eof

For ease of notation, let $\widetilde{p}^{\,k+1}_j = p^{k+1}r_jv_j - c_j$ denote the profit of product $j \in P$ in the interval $[p^k,p^{k+1}] \in \mathcal{I}(\gexp)$. As previously mentioned, computing the parametric bound $\ovzpk$ requires sorting the products by non-increasing utility-to-space consumption ratios $\widetilde{p}^{\,k +1}_j/ v_j$ and filling up the knapsack starting from a product with the largest ratio. For each interval $[p^k, p^{k+1}]$, we can thus keep track of the critical product $s(k)$, i.e., the largest index such that $\sum_{j=1}^{s(k)-1} v_j \leq 1/p^k -1$. We can then state the following variable-fixing criterion.

\begin{proposition}
For a given $\rho$, each product $j \in P$ such that, for all intervals $[p^k,p^{k+1}] \in \mathcal{I}(\gexp)$, the following two conditions hold $$(i) \;\; \widetilde{p}^{\,k +1}_j/ v_j < \widetilde{p}^{\,k+1}_{s(k)}/ v_{s(k)} \quad \text{and} \quad (ii) \;\; \ovzpk + \widetilde{p}^{\,k+1}_j -  v_{j}\widetilde{p}^{\,k+1}_{s(k)}/ v_{s(k)} < lb(\gexp)$$
cannot be part of any \reftwo{optimal \probname assortment}.
\end{proposition}
\textit{Proof.} Let $\ov{x}$ denote the optimal solution of problem (\ref{FTparub}) with objective value $\ovzpk$, for a given $[p^k,p^{k+1}] \in \mathcal{I}(\gexp)$. We observe that $\ov{x}_j=0$ for each product $j\in P$ satisfying the first condition (since its utility-to-space ratio is smaller than the one of the critical product $s(k)$) and that, {for such products}, the quantity $\widetilde{p}^{\,k+1}_j -  v_{j}\widetilde{p}^{\,k+1}_{s(k)}/ v_{s(k)}$ is a lower bound of the profit decrease obtained by setting $\ov{x}_j=1$ \citep[see, e.g.,][Section 2.2.3]{KPbooktothMartello}. {The following inequality thus holds:
\begin{equation*}
    \ovzpk \geq \ovzpk + \widetilde{p}^{\,k+1}_j -  v_{j}\widetilde{p}^{\,k+1}_{s(k)}/ v_{s(k)} \geq \ovzpk_{|x_j=1},
\end{equation*}
where $\ovzpk_{|x_j=1}$ is the solution of the Linear Knapsack relaxation \eqref{FTparub} obtained after fixing $x_j$ to 1, given the interval $[p^k,p^{k+1}] \in \mathcal{I}(\gexp)$. The proposition follows by noticing that condition (ii) implies $lb(\gexp) \geq  \ovzpk_{|x_j =1} \text{ for all } [p^k,p^{k+1}] \in \mathcal{I}(\gexp)$. In words,}  introducing product $j$ in any of the solutions corresponding to $[p^k,p^{k+1}] \in \mathcal{I}(\gexp)$ results in a profit lower than the best-known feasible solution. Hence, product $j$ cannot be part of any \reftwo{optimal \probname assortment}. \eof

\subsection{Finding an optimal \probname assortment} \label{subsec:gap}
The bounding procedure described in \cref{subsec:seqbound} returns two important values, namely $\pmin(\gexp)$ and $\pmax(\gexp)$, for $\rho = \rho^\ell$, that identify a crucial feature of any \reftwo{optimal \probname assortment}: the corresponding no-purchase probability lies in the range $[\pmin(\gexp),\pmax(\gexp)]$. To identify an \reftwo{optimal \probname assortment}, we solve problem \eqref{prob:milp}, with a general-purpose solver, with the addition of the constraints 
\begin{equation}
\pmin(\gexp) \leq u_0 \leq \pmax(\gexp),
\label{addCon}
\end{equation}
which limit the search for an \reftwo{optimal \probname assortment} to assortments having no-purchase probability within the range $[\pmin(\gexp),\pmax(\gexp)]$. 

\reftwo{Moreover, because of the variable fixing rules described in Section \ref{subsec:varfix}, the $x_j$ variables corresponding to the items that can be ruled out are set equal to 0.}

In the following, we refer to problem \eqref{prob:milp} plus constraints \eqref{addCon} \reftwo{and constraint $x_j = 0$ for each item $j$ that can be ruled out} as \final. 

\reftwo{An optimal \probname assortment is found by solving \final\ with a general-purpose MILP solver. We could have also investigated how to find an optimal \probname assortment by exploiting the MICQ model \eqref{SAK}. 
We did not explore this alternative because the computing time to solve \final\ is, on average, about a tenth of a second even on large instances as will be shown in the computational experiments summarized in Section \ref{sec:results}.} 

\vspace{0.5cm}

{We conclude this section by summarizing the three steps of our exact solution method (hereafter referred to as LLRS\footnote{{From the initials of the authors' family names.}}) i.e., enhanced bounding procedure, variable fixing and search for the optimal assortment, respectively, which we report in the Pseudocode \ref{pseudo}.  

\begin{algorithm}\label{pseudo}
{
\caption{Pseudocode of LLRS}
\SetKwInOut{Input}{Input}
\SetKwInOut{Output}{Output}
\Input{\rhof, \rhol}
\Output{Optimal Solution}
$\rho \gets \rho^f$\;
\texttt{Ruled\_out}$=\emptyset$\;
$p_{min} = 1/v_0$\;
$p_{max} = 1/(v_0 + \pmb{1}^T\pmb{v})$\;
\tcp {Enhanced Bounding Procedure}
\While{$\rho\geq$\rhol}{
    $z(\gexp)\gets 0$\;
    $lb(\gexp)\gets 0$\;
    $\mathcal{I}(\gexp) = \emptyset$\;
    \tcp{Compute upper and lower bounds}
    \ForEach {$k \in \{1,\dots,K(\rho)+1\}$}{
        $ ub \gets \ovzpk$ \tcp*{Linear Knapsack \eqref{FTparub} with optimal solution $\ov{x}$}
        $ lb \gets \unzpk$ \tcp*{Obtained removing critical item $s(k)$ from $\ov{x}$}
        $z(\gexp) \gets max\{z(\gexp), ub\}$ \;
        $lb(\gexp) \gets max\{lb(\gexp), lb\}$ \;
    }
    \tcp{Prune intervals}
    \ForEach {$k \in \{1,\dots,K(\rho)+1\}$}{
        \If{$\ovzpk \geq lb(\gexp)$}{
            $I(\gexp) = \mathcal{I}(\gexp)\cup \{(p^k,p^{k+1})\}$\;
        }
    }
    $p_{min} \gets p_{min}(\gexp)$\;
    $p_{max} \gets p_{max}(\gexp)$\;
    \tcp{Reductions}
    \ForEach {$j \in P\setminus$ \textup{\texttt{Ruled\_out}}}{
        \If{Proposition 1 or Proposition 2 apply}{
            \texttt{Ruled\_out} = \texttt{Ruled\_out}$\cup \{j\}$
        }
    }
    $\rho \gets \rho / 10$\;
}
    \SetKwFunction{FMain}{MILP+}
      \SetKwProg{Fn}{Execute}{:}{}
      \Fn{\FMain{$p_{min}$, $p_{max}$, \textup{\texttt{Ruled\_out}}}}{
      }
    }
\end{algorithm}
}

\section{Constrained \probname}\label{sec:constrained}
In real-life applications, there could be additional constraints that should be taken into account when selecting an assortment. The two most common examples of such constraints are \textit{cardinality constraints} and \textit{space/resource constraints} \citep[see, e.g.,][]{Feldman2015,Feldman2019,Kunnumkal2019}. Here, we show how the exact method described in \cref{sec:framework} can be adjusted to handle a cardinality constraint. A space/resource capacity constraint can be handled similarly.
A cardinality constraint implies that no more than a given number of $\kappa$ products can be included in any assortment, i.e., 
\begin{equation}
\sum_{j \in P} x_j \leq \kappa ~ . \label{cardcon}
\end{equation}

The following three modifications are made to the exact method to consider constraint \eqref{cardcon}:
\begin{enumerate}
    \item The presence of constraint \eqref{cardcon} requires solving the linear relaxation of a two-constrained binary knapsack problem \citep[see][]{martello1997,martello2003} when computing the parametric \refone{upper} bound $\ovzpk$. \reftwo{\cite{martello1997} prove that $\ovzpk$ can be computed in polynomial time $O(n^2)$ by relaxing constraint \eqref{cardcon} through a non-negative Lagrangian multiplier $\lambda \in \R$; however, we propose a different approach and compute an (experimentally tight) upper bound to $\ovzpk$ by considering a limited set of non-decreasing values of $\lambda$, which does not necessarily include the optimal value $\lambda^\ast$, when relaxing constraint \eqref{cardcon}. We note, however, that the results obtained in preliminary computational experiments suggest that optimal Lagrangian multipliers $\lambda^\ast$ seem to be non-decreasing in $k$. Thus, we exploit the similarity of the two-constrained binary knapsack problems whose linear relaxations need to be solved to compute $z(\gexp)$ (i.e., $K(\rho)$ problems in total) by calculating an upper bound $\ov{z}(p^k,p^{k+1},\lambda)$ to $\ovzpk$ using the following procedure.}
    \reftwo{For each $k \in \{1,\ldots,K(\rho)\}$, $\ovzpk$ is set equal to the minimum over a limited set of values $\ov{z}(p^k,p^{k+1},\lambda)$ with $\lambda \in \R$, where $\ov{z}(p^k,p^{k+1},\lambda)$ is defined as}
\begin{equation}
\reftwo{\ov{z}(p^k,p^{k+1},\lambda) = \max_{\gb{x} \in \BL^n} \Bigg\{ \sum_{j \in P} \left( p^{k+1} r_j v_j - c_j - \lambda \right)x_j + \lambda \kappa \, | \, \vx \leq \frac{1}{p^{k}} - 1 \Bigg\}.}
\end{equation}
For a given grid density $\rho$, we initialize $\lambda$ to 0. Starting from the interval corresponding to $p^k \leq \pmax(\gexp) \leq p^{k+1}$, we compute a ``good'' Lagrangian multiplier (and upper bound) iteratively by using, for iteration $t$, $\lambda_t=\lambda_{t-1} + \delta$, with $\delta$ small enough. We stop the procedure at iteration $\ov{t}$ if the upper bound starts deteriorating, i.e., if $\ov{z}(p^k,p^{k+1},\lambda_{\ov{t}}) > \ov{z}(p^k,p^{k+1},\lambda_{\ov{t}-1})$. We then use the close-to-optimal multiplier $\lambda_{\ov{t}}$ to initialize the procedure for the next interval, thus avoiding starting from scratch, i.e., by assuming $\lambda = 0$. The same approach is used for all other intervals, i.e., until $p^k \leq \pmin(\gexp) \leq p^{k+1}$. Experimentally, we found this procedure to recover good multipliers and upper bounds in a small number of iterations.
    \item Products are added to assortment $P'$ as long as the sum of the preferences of the selected products does not exceed $\frac{1}{\unp} - 1$ \textit{and} if the total number of products $P'$ is not greater than $\kappa$ when computing $\unzpk$ (see \cref{subsec:seqbound}).
    \item Constraint \eqref{cardcon} is added to \final\ when computing an \reftwo{optimal \probname assortment} (see \cref{subsec:gap}).
\end{enumerate}

\section{Computational results\label{sec:results}}
In this section, we first describe the instances used to test our exact method (see \cref{subsec:inst}). Then, we report the results achieved by our exact method on the \probname (see \cref{subsec:resunc}) and the cardinality-constrained \probname (see \cref{subsec:rescon}). We also compare the results of our exact method with those of the two formulations from the literature described in \cref{sec:prob}. \blue{All the code and instances used to run the experiments are publicly available \footnote{\texttt{https://github.com/ds4dm/LLRS}}.}

\subsection{Test instances} \label{subsec:inst}
To test our exact method, we generate a set of 800 instances as described by \citet{Kunnumkal2019}. Each instance has $n \in \{100, 200, 500, 1000\}$ products. The preference $v_j$ of product $j \in P$ is computed as $v_j = w_j / \sum_{k=1}^n w_k$, where $w_j$ is uniformly distributed in the interval $(0,1]$. The no-purchase probability is set equal to $v_0 = \frac{\varPhi}{1-\varPhi}\sum_{j \in P} v_j$, where $\varPhi$ is a parameter in the set $\varPhi \in \{0.25,0.75\}$, meaning that the no-purchase probability is either 25\% or 75\% {when all products are included in the assortment}. The revenue $r_j$ of product $j \in P$ is sampled from the uniform distribution $[0,2000]$, and the product cost $c_j$ is sampled from the uniform distribution $[0,\gamma r_j v_j / (v_0 + v_j)]$, where $\gamma$ is another parameter in the set $\gamma \in \{0.5, 1.0\}$. Therefore, we have 16 combinations of instances, $\comb$. For each combination, we generate 50 test instances.

\subsection{Computational results on the \probname} \label{subsec:resunc}
Table \ref{tab:summ} summarizes the computational results achieved with our {LLRS} 
on the 800 test instances and compares its performance with the performance of MILP \eqref{prob:milp} and MICQ \eqref{SAK} both solved with Cplex 20.1. All experiments are conduced on a single core of a machine with 500GB-RAM and an Intel(R) Xeon(R)Gold 6142 with 2.60GHz CPU. A time limit of ten minutes is imposed on each experiment.


{
Our implementation has been coded in C++ and compiled with \texttt{gcc 4.8.5}; Cplex 20.1 is used to solve problem \final. Our exact method has two parameters: $\rho^f$ and $\rho^\ell$ that we set as $\rho^f = \num{1e-2}$ and $\rho^\ell = \num{1e-7}$, respectively. In this regard, we conducted a set of experiments 
whose results assessed the robustness of LLRS with respect to these parameters and helped identifying simple rules for setting them. In particular, one may want to set \rhof relatively large, since (i) coarse grids help pruning more refined ones, and (ii) even if large values of \rhof may provide loose bounds, they are extremely fast to compute. In practice, values of \rhof such that the number of intervals (i.e., linear knapsack relaxations) to solve at the first iterations is about one hundred or smaller performed well. Concerning parameter \rhol, it should be set relatively small as smaller values decrease the average computing time of \final. In this regards, we observed that our bounding procedure provides a significant computational boost for $\rho^\ell\le \num{1e-4}$.} \blue{We refer the interested reader to Appendix \ref{ap:sensitivity}, where we report and discuss the results of this set of experiments.}

For each of the three methods and each combination of test instances, Table \ref{tab:summ} reports the number of instances solved to optimality (\textsf{opt}), the average final gap (\textsf{gap}) in percentage between the best \refone{lower} and \refone{upper} bounds found by the corresponding method, the average computing time (\textsf{\cpuavg}), and the maximum computing time (\textsf{\cpumax}).

\begin{table}[ht]
\begin{small}
\begin{center}
\sffamily
\caption{Summary of the computational results of \our\ on the \probname and comparison with MILP \eqref{prob:milp} and MICQ \eqref{SAK}} \label{tab:summ}
\setlength{\tabcolsep}{5.5pt}
\renewcommand{\arraystretch}{1.1}
\begin{tabular}{crrrrrrrrrrr}
\toprule
 & \multicolumn{4}{c}{MILP \eqref{prob:milp}} & \multicolumn{4}{c}{MICQ \eqref{SAK}} & \multicolumn{3}{c}{LLRS} \\
\cmidrule(r){2-5} \cmidrule(r){6-9} \cmidrule{10-12}
 $\comb$ & opt & gap & \cpuavg & \cpumax & opt & gap & \cpuavg & \cpumax    & opt & \cpuavg & \cpumax \\
\midrule
(100, 0.25, 0.5)  & 48  &  0.08 & 112.46 & 600.00 & 14  &  2.20 & 493.46 & 600.00 & 50  & 0.02 & 0.10 \\
(100, 0.25, 1.0)  & 50  &  0.00 &   1.40 &   4.90 & 43  &  0.35 & 185.84 & 600.00 & 50  & 0.01 & 0.08 \\
(100, 0.75, 0.5)  & 50  &  0.00 &   0.20 &   0.32 & 50  &  0.00 &   2.28 &  10.46 & 50  & 0.00 & 0.02 \\
(100, 0.75, 1.0)  & 50  &  0.00 &   0.29 &   0.80 & 48  &  0.03 &  38.22 & 600.00 & 50  & 0.00 & 0.02 \\
\tabsp
(200, 0.25, 0.5)  & 0   & 11.86 & 600.00 & 600.00 &  0  & 11.02 & 600.00 & 600.00 & 50  & 0.02 & 0.09 \\
(200, 0.25, 1.0)  & 10  &  3.83 & 546.22 & 600.00 &  0  &  9.39 & 600.00 & 600.00 & 50  & 0.01 & 0.06 \\
(200, 0.75, 0.5)  & 50  &  0.00 &   0.33 &   0.53 & 42  &  0.04 & 161.54 & 600.00 & 50  & 0.01 & 0.04 \\
(200, 0.75, 1.0)  & 47  &  0.01 &  82.21 & 600.00 &  0  &  1.58 & 600.00 & 600.00 & 50  & 0.01 & 0.04 \\
\tabsp
(500, 0.25, 0.5)  & 0   & 24.84 & 600.00 & 600.00 &  0  & 25.20 & 600.00 & 600.00 & 50  & 0.06 & 0.49 \\
(500, 0.25, 1.0)  & 0   & 20.69 & 600.00 & 600.00 &  0  & 26.38 & 600.00 & 600.00 & 50  & 0.04 & 0.19 \\
(500, 0.75, 0.5)  & 50  &  0.00 &  18.25 & 196.54 &  0  &  0.54 & 600.00 & 600.00 & 50  & 0.05 & 0.10 \\
(500, 0.75, 1.0)  & 0   &  1.34 & 600.00 & 600.00 &  0  &  3.66 & 600.00 & 600.00 & 50  & 0.04 & 0.09 \\
\tabsp
(1000, 0.25, 0.5) & 0   & 28.39 & 600.00 & 600.00 &  0  & 25.74 & 600.00 & 600.00 & 50  & 0.29 & 1.78 \\
(1000, 0.25, 1.0) & 0   & 26.69 & 600.00 & 600.00 &  0  & 39.00 & 600.00 & 600.00 & 50  & 0.15 & 1.06 \\
(1000, 0.75, 0.5) & 0   &  0.07 & 600.00 & 600.00 &  0  &  0.92 & 600.00 & 600.00 & 50  & 0.16 & 0.24 \\
(1000, 0.75, 1.0) & 0   &  2.15 & 600.00 & 600.00 &  0  &  4.33 & 600.00 & 600.00 & 50  & 0.12 & 0.22 \\
\midrule
Sum & 355 & \multicolumn{3}{c}{} & 197 & \multicolumn{3}{c}{} & 800 \\	
Avg & & 7.50 & 347.58 & & & 9.40 & 467.58 & & & 0.06 & \\
Max & & 28.39 & & 600.00 & & 39.00 & & 600.00 & & & 1.78 \\
\bottomrule
\end{tabular}
\end{center}
\end{small}
\end{table}

Table \ref{tab:summ} shows that \our\ outperforms the two formulations from the literature in terms of both number of instances solved to optimality and computing time. Indeed, \our\ can solve all 800 instances whereas the MILP \eqref{prob:milp} and the MICQ \eqref{SAK} can solve 355 and 197 instances only, respectively. The average computing time of \our\ is at least four/five orders of magnitude lower than the computing time of the other two methods. We can observe that the average computing time of \our\ increases with the number of products and that instances featuring lower product costs (i.e., $\gamma = 0.5$) or lower no-purchase probability (i.e., $\varPhi = 0.25$) are more difficult to solve than instances with higher product costs (i.e., $\gamma = 1.0$) or higher no-purchase probability (i.e., $\varPhi = 0.75$). It is interesting to observe that the average computing time of \our\ on the large instances with 1000 products is about 0.2 seconds and the maximum computing time is just 1.78 seconds.

Table \ref{tab:addInfo} reports additional information on the performance of \our\ for each combination of test instances. 
For the bounding procedure described in \cref{subsec:seqbound}, Table \ref{tab:addInfo} reports the final average \refone{upper} gap (\textsf{gap$_{\text{dual}}$}) in percentage (i.e., the gap between $z^\ast$ and $z(\gexp)$ with $\rho = \rho^\ell$), the final average \refone{lower} gap (\textsf{gap$_{\text{prim}}$}) in percentage (i.e., the gap between $z^\ast$ and $lb(\gexp)$ with $\rho = \rho^\ell$), the number of times $lb(\gexp)$ corresponds to the optimal solution (\textsf{opt$_{\text{prim}}$}), and the average computing time (\textsf{\cpuavg}). For the variable fixing described in \cref{subsec:varfix}, \textsf{\#out} indicates the average number of products ruled out and \textsf{\%out} is the number of products ruled out in percentage. Finally, for the resolution of \final\ (see \cref{subsec:gap}), \refone{\textsf{gap$_{\text{LP}}$} is the average gap in percentage between the lower bound provided by its linear relaxation and $z^\ast$, and} \textsf{\cpuavg} is the average computing time\refone{; under label \textsf{gap$_{\text{LP}}^{\text{MILP}\eqref{prob:milp}}$}, we report the average gap in percentage between the linear relaxation of MILP \eqref{prob:milp} and $z^\ast$ to have a comparison with \textsf{gap$_{\text{LP}}$}.}

\begin{table}[ht]
\begin{small}
\begin{center}
\sffamily
\caption{Additional information on the performance of \our\ on the \probname} \label{tab:addInfo}
\setlength{\tabcolsep}{5.5pt}
\renewcommand{\arraystretch}{1.1}
\begin{tabular}{crrrrrrrrrr}
\toprule
 & \multicolumn{4}{c}{Bounding Procedure} & \multicolumn{2}{c}{Var Fixing} & \refone{MILP} & \multicolumn{2}{c}{\final} \\
\cmidrule(r){2-5} \cmidrule(r){6-7} \cmidrule(r){8-8} \cmidrule(r){9-10}
$\comb$ & gap$_{\text{dual}}$ & gap$_{\text{prim}}$ & opt$_{\text{prim}}$ & \cpuavg & \#out & \%out & \refone{gap$_{\text{LP}}$} & \refone{gap$_{\text{LP}}$} & \cpuavg \\
\midrule
(100, 0.25, 0.5) & 0.0031 & 0.0000 & 49 & 0.02 &  62.2 & 62.2 & 31.3 & 18.2 & 0.00 \\
(100, 0.25, 1.0) & 0.0047 & 0.0001 & 49 & 0.01 &  72.4 & 72.4 & 31.2 & 18.1 & 0.00 \\
(100, 0.75, 0.5) & 0.0000 & 0.0000 & 50 & 0.00 &   0.9 &  0.9 & 10.0 &  9.5 & 0.00 \\
(100, 0.75, 1.0) & 0.0002 & 0.0000 & 50 & 0.00 &  14.7 & 14.7 & 16.0 & 12.1 & 0.00 \\
\tabsp
(200, 0.25, 0.5) & 0.0007 & 0.0001 & 48 & 0.02 & 113.6 & 56.8 & 31.4 & 19.0 & 0.00 \\
(200, 0.25, 1.0) & 0.0008 & 0.0000 & 49 & 0.01 & 137.9 & 68.9 & 32.2 & 19.2 & 0.00 \\
(200, 0.75, 0.5) & 0.0000 & 0.0000 & 50 & 0.00 &   0.0 &  0.0 & 10.0 &  9.5 & 0.01 \\
(200, 0.75, 1.0) & 0.0001 & 0.0000 & 50 & 0.00 &  26.5 & 13.3 & 16.4 & 12.5 & 0.01 \\
\tabsp
(500, 0.25, 0.5) & 0.0001 & 0.0000 & 49 & 0.03 & 221.6 & 44.3 & 32.0 & 19.7 & 0.03 \\
(500, 0.25, 1.0) & 0.0002 & 0.0000 & 47 & 0.02 & 304.1 & 60.8 & 31.3 & 19.4 & 0.02 \\
(500, 0.75, 0.5) & 0.0000 & 0.0000 & 50 & 0.01 &   0.0 &  0.0 & 10.1 &  9.6 & 0.04 \\
(500, 0.75, 1.0) & 0.0000 & 0.0000 & 50 & 0.01 &  65.5 & 13.1 & 16.3 & 12.4 & 0.03 \\
\tabsp
(1000, 0.25, 0.5) & 0.0000 & 0.0000 & 48 & 0.12 & 234.6 & 23.5 & 31.8 & 20.2 & 0.17 \\
(1000, 0.25, 1.0) & 0.0001 & 0.0000 & 49 & 0.06 & 483.3 & 48.3 & 31.8 & 20.3 & 0.09 \\
(1000, 0.75, 0.5) & 0.0000 & 0.0000 & 50 & 0.04 &   0.0 &  0.0 & 10.1 &  9.6 & 0.12 \\
(1000, 0.75, 1.0) & 0.0000 & 0.0000 & 50 & 0.04 & 135.0 & 13.5 & 16.1 & 12.4 & 0.08 \\
\midrule
Avg & 0.0006 & 0.0000 &  & 0.03 & 117.0 & 30.8 & \refone{22.4} & \refone{15.1} & 0.04 \\
\bottomrule
\end{tabular}
\end{center}
\end{small}
\end{table}

Table \ref{tab:addInfo} shows that both the \refone{lower} and \refone{upper} bounds computed by the bounding procedure are very tight and the \refone{lower} bound $lb(\gexp)$ is most of the times an \reftwo{optimal \probname assortment} (in all but 12 instances). We can also observe that roughly half of the computing time of \our\ is spent on the bounding procedure and the other half on solving \final. 

\reftwo{Table \ref{tab:addInfo} also indicates that the proposed variable fixing rules allow to rule out about 31\% of the products on average. We do not report detailed results on solving the \probname with \our\ without applying the variable fixing rules. However, we noticed that the computing time was, on average, around 30\% higher than when the variable fixing rules were applied.}

\refone{Finally, Table \ref{tab:addInfo} shows that the linear relaxation of \final\ is, on average, much tighter than the linear relaxation of MILP \eqref{prob:milp}: indeed, the average percentage gap is 15.1\% versus 22.4\%. This is due to the addition of constraints \eqref{addCon}.}

We conduct another set of experiments to show the importance of applying the bounding procedure described in \cref{subsec:seqbound}. In particular, we compare the results achieved by setting $\rho^f$ either to $\num{1E-02}$ or to $\num{1E-07}$. When $\rho^f = \num{1E-02}$, we have the results already reported in Table \ref{tab:summ}. When $\rho^f = \num{1E-07}$, the bounding procedure computes a single \refone{upper} bound $z(\gexp)$, with $\rho = \num{1E-07}$, and cannot benefit from the bounds achieved with higher values of $\rho$; in other words, this approach corresponds to the one proposed by \citet{Feldman2015}. Table \ref{tab:vsFT} reports the average and maximum number of intervals for which $\ovzpk$ must be computed (\textsf{\#\intavg}, \textsf{\#\intmax} - notice that \textsf{\#\intavg} is equal to \textsf{\#\intmax} when $\rho^f = \num{1E-07}$), the average and maximum computing time to solve the problem (\textsf{\cpuavg}, \textsf{\cpumax}), the average and maximum number of intervals, in percentage, for which $\ovzpk$ is computed when $\rho^f = \num{1E-02}$ compared to when $\rho^f = \num{1E-07}$.

\begin{table}[ht]
\begin{small}
\begin{center}
\sffamily
\caption{Comparison with the results achieved by \our\ by setting $\rho^f = \num{1E-07}$} \label{tab:vsFT}
\setlength{\tabcolsep}{5.5pt}
\renewcommand{\arraystretch}{1.1}
\begin{tabular}{crrrrrrrrr}
\toprule
 & \multicolumn{3}{c}{\our\ ($\rho^f = 1\mathrm{E}-7$)} & \multicolumn{6}{c}{\our\ ($\rho^f = 1\mathrm{E}-2$)} \\
\cmidrule(r){2-4} \cmidrule{5-10}
$\comb$ & \#\intavg & \cpuavg & \cpumax & \#\intavg & \%\intavg & \#\intmax & \%\intmax & \cpuavg & \cpumax \\
\midrule
(100, 0.25, 0.5) & 13\,862\,947 &  7.04 &  7.19 & 62\,251 & 0.4 & 299\,206 & 2.2 & 0.02 & 0.10 \\
(100, 0.25, 1.0) & 13\,862\,947 &  7.36 &  7.52 & 57\,430 & 0.4 & 348\,920 & 2.5 & 0.01 & 0.08 \\
(100, 0.75, 0.5) &  2\,876\,823 &  1.44 &  1.47 &  2\,906 & 0.1 &  30\,049 & 1.0 & 0.00 & 0.02 \\
(100, 0.75, 1.0) &  2\,876\,823 &  1.46 &  1.49 &  4\,657 & 0.2 &  40\,371 & 1.4 & 0.00 & 0.02 \\
\tabsp
(200, 0.25, 0.5) & 13\,862\,947 & 13.24 & 13.58 & 37\,499 & 0.3 & 142\,049 & 1.0 & 0.02 & 0.09 \\
(200, 0.25, 1.0) & 13\,862\,947 & 13.86 & 14.18 & 32\,151 & 0.2 & 154\,331 & 1.1 & 0.01 & 0.06 \\
(200, 0.75, 0.5) &  2\,876\,823 &  2.71 &  2.80 &  4\,421 & 0.2 &  29\,177 & 1.0 & 0.01 & 0.04 \\
(200, 0.75, 1.0) &  2\,876\,823 &  2.73 &  2.78 &  5\,871 & 0.2 &  30\,551 & 1.1 & 0.01 & 0.04 \\
\tabsp
(500, 0.25, 0.5) & 13\,862\,947 & 40.00 & 42.14 & 23\,993 & 0.2 &  64\,459 & 0.5 & 0.06 & 0.49 \\
(500, 0.25, 1.0) & 13\,862\,947 & 42.62 & 43.88 & 22\,651 & 0.2 &  72\,100 & 0.5 & 0.04 & 0.19 \\
(500, 0.75, 0.5) &  2\,876\,823 &  7.40 &  7.51 &  5\,028 & 0.2 &  21\,586 & 0.8 & 0.05 & 0.10 \\
(500, 0.75, 1.0) &  2\,876\,823 &  7.79 &  8.16 &  5\,649 & 0.2 &  15\,801 & 0.5 & 0.04 & 0.09 \\
\tabsp
(1000, 0.25, 0.5) & 13\,862\,947 & 91.33 & 93.87 & 22\,701 & 0.2 &  42\,799 & 0.3 & 0.29 & 1.78 \\
(1000, 0.25, 1.0) & 13\,862\,947 & 98.27 &100.98 & 19\,037 & 0.1 &  46\,461 & 0.3 & 0.15 & 1.06 \\
(1000, 0.75, 0.5) &  2\,876\,823 & 16.79 & 17.08 &  5\,749 & 0.2 &  16\,638 & 0.6 & 0.16 & 0.24 \\
(1000, 0.75, 1.0) &  2\,876\,823 & 17.93 & 18.29 &  6\,076 & 0.2 &  16\,971 & 0.6 & 0.12 & 0.22 \\
\midrule
Avg & 8\,369\,885 & 23.25 & & 19\,879 & 0.2 & & & 0.06 \\	
Max & & & 100.98 & & & 348\,920 & 2.5 & & 1.78 \\
\bottomrule
\end{tabular}
\end{center}
\end{small}
\end{table}

Table \ref{tab:vsFT} shows that applying the bounding procedure as described in \cref{subsec:seqbound} allows to significantly reduce the number of times the parametric bound must be computed: on average by 99.8\%. This translates into much lower computing times, which decrease by two/three orders of magnitude on average. We can also observe that \textsf{\#\intavg} gradually decreases when the number of products increases and, not surprisingly, when the no-purchase option is higher. 

\refone{The computational experiments summarized in Tables \ref{tab:summ}-\ref{tab:vsFT} have been conducted on a set of instances generated as in \cite{Kunnumkal2019}, where the no-purchase probability ($\varPhi$) is either 25\% or 75\%. The results show that, when $\varPhi$ equals $0.25$, the average computing time to solve the \probname\ to optimality with \our\ as well as the average dual gaps of the bounding procedure (\textsf{gap$_{\text{dual}}$}) and the linear relaxation of \final\ (\textsf{gap$_{\text{LP}}$}) are all higher than the case where $\varPhi$ equals $0.75$. This behavior may suggest that if $\varPhi$ is even lower, \our\ may not be able to solve \probname\ instances of the same size within the same amount of computing time. Hence, we conducted another set of experiments on the same instances previously considered but setting $\varPhi$ equal to 0.1. The results are summarized in Table \ref{tab:summAdd} that reports the same information as Table \ref{tab:summ}.}

\begin{table}[ht]
\begin{small}
\begin{center}
\sffamily
\caption{\refone{Summary of the computational results of \our\ on the \probname with $\varPhi = 0.1$ and comparison with MILP \eqref{prob:milp} and MICQ \eqref{SAK}}} \label{tab:summAdd}
\setlength{\tabcolsep}{5.5pt}
\renewcommand{\arraystretch}{1.1}
\begin{tabular}{crrrrrrrrrrr}
\toprule
 & \multicolumn{4}{c}{MILP \eqref{prob:milp}} & \multicolumn{4}{c}{MICQ \eqref{SAK}} & \multicolumn{3}{c}{LLRS} \\
\cmidrule(r){2-5} \cmidrule(r){6-9} \cmidrule{10-12}
 $\comb$ & opt & gap & \cpuavg & \cpumax & opt & gap & \cpuavg & \cpumax    & opt & \cpuavg & \cpumax \\
\midrule
(100, 0.1, 0.5)  & 50  &  0.00 &   4.52 & 100.34 & 49  &  0.28 &  32.45 & 600.00 & 50  & 0.03 & 0.12 \\
(100, 0.1, 1.0)  & 50  &  0.00 &   0.38 &   4.85 & 50  &  0.00 &   9.23 &  83.30 & 50  & 0.03 & 0.17 \\
\tabsp
(200, 0.1, 0.5)  & 10  &  4.81 & 530.48 & 600.00 &  1  &  7.24 & 592.75 & 600.00 & 50  & 0.03 & 0.11 \\
(200, 0.1, 1.0)  & 46  &  0.25 & 110.21 & 600.00 & 10  &  7.60 & 549.18 & 600.00 & 50  & 0.03 & 0.11 \\
\tabsp
(500, 0.1, 0.5)  & 0   & 21.46 & 600.00 & 600.00 &  0  & 21.56 & 600.00 & 600.00 & 50  & 0.05 & 0.27 \\
(500, 0.1, 1.0)  & 0   & 15.97 & 600.00 & 600.00 &  0  & 23.85 & 600.00 & 600.00 & 50  & 0.03 & 0.19 \\
\tabsp
(1000, 0.1, 0.5) & 0   & 27.41 & 600.00 & 600.00 &  0  & 24.41 & 600.00 & 600.00 & 50  & 0.14 & 0.91 \\
(1000, 0.1, 1.0) & 0   & 26.40 & 600.00 & 600.00 &  0  & 30.88 & 600.00 & 600.00 & 50  & 0.05 & 0.27 \\
\midrule
Sum & 156 & \multicolumn{3}{c}{} & 110 & \multicolumn{3}{c}{} & 400 \\	
Avg & & 12.04 & 380.68 & & & 14.48 & 447.96 & & & 0.05 & \\
Max & & 27.41 & & 600.00 & & 30.88 & & 600.00 & & & 0.91 \\
\bottomrule
\end{tabular}
\end{center}
\end{small}
\end{table}

\refone{Table \ref{tab:summAdd} shows that, even when $\varPhi$ equals 0.1, \our\ can still solve all instances to optimality with an average computing time of 0.05 seconds and a maximum computing time of 0.91 seconds. \our\ still outperforms both MILP \eqref{prob:milp} and the MICQ \eqref{SAK}, which can solve 156 and 110 instances, respectively. 
}

\subsection{Computational results on the cardinality-constrained \probname} \label{subsec:rescon}
In this section, we report on the performance of \our\ on the cardinality-constrained \probname. We use the same 800 instances as in the computational study on the \probname and set the maximum cardinality of the assortments equal to half of the number of products, i.e., $\kappa = n / 2$, as done in \citet{Kunnumkal2019}. For \our, parameter $\delta$ is set equal to \num{1E-5}.

Table \ref{tab:summConstr} summarizes the results achieved on these 800 instances by \our\ in the same format of Table \ref{tab:summ} and compares them to the results of MILP \eqref{prob:milp} and MICQ \eqref{SAK} both solved with Cplex 20.1.

Table \ref{tab:summConstr} shows that \our\ outperforms the two formulations from the literature as it can solve all of the 800 instances to optimality, compared to 302 instances solved by the MILP and 354 instances solved by the MICQ, and is orders of magnitude faster than these two formulations. As observed on the \probname, the computing time of \our\ increases with the number of products. Furthermore, instances featuring lower product costs (i.e., $\gamma = 0.5$) or higher no-purchase probability (i.e., $\varPhi = 0.75$) take more time to solve to optimality. It is interesting to observe that the average computing time of \our\ on large instances with 1000 products is approximately 0.67 seconds, and the maximum computing time over all instances is just 2.27 seconds.

\begin{table}[!htbp]
\begin{small}
\begin{center}
\sffamily
\caption{Summary of the computational results of \our\ on the cardinality-constrained \probname and comparison with MILP \eqref{prob:milp} and MICQ \eqref{SAK}} \label{tab:summConstr}
\setlength{\tabcolsep}{5.5pt}
\renewcommand{\arraystretch}{1.1}
\begin{tabular}{crrrrrrrrrrr}
\toprule
 & \multicolumn{4}{c}{MILP \eqref{prob:milp}} & \multicolumn{4}{c}{MICQ \eqref{SAK}} & \multicolumn{3}{c}{LLRS} \\
\cmidrule(r){2-5} \cmidrule(r){6-9} \cmidrule{10-12}
$\comb$ & opt & gap & \cpuavg & \cpumax & opt & gap & \cpuavg & \cpumax    & opt & \cpuavg & \cpumax \\
\midrule
(100, 0.25, 0.5) & 48 &  0.08 & 115.41 & 600.00 & 50 &  0.00 &  27.15 & 250.47 & 50 & 0.02 & 0.09 \\
(100, 0.25, 1.0) & 50 &  0.00 &   1.50 &   5.03 & 49 &  0.05 &  62.04 & 600.00 & 50 & 0.01 & 0.08 \\
(100, 0.75, 0.5) & 50 &  0.00 &   0.39 &   1.25 & 50 &  0.00 &   0.39 &   0.87 & 50 & 0.08 & 0.12 \\
(100, 0.75, 1.0) & 50 &  0.00 &   0.40 &   0.84 & 50 &  0.00 &   0.94 &   2.92 & 50 & 0.02 & 0.05 \\
\tabsp
(200, 0.25, 0.5) & 0  & 11.95 & 600.00 & 600.00 &  2 &  1.68 & 586.84 & 600.00 & 50 & 0.03 & 0.10 \\
(200, 0.25, 1.0) & 9  &  3.82 & 550.88 & 600.00 &  1 &  4.81 & 595.49 & 600.00 & 50 & 0.01 & 0.06 \\
(200, 0.75, 0.5) & 48 &  0.00 & 115.14 & 600.00 & 50 &  0.00 &   1.28 &  10.72 & 50 & 0.17 & 0.24 \\
(200, 0.75, 1.0) & 47 &  0.02 &  81.94 & 600.00 & 50 &  0.00 &  25.76 & 299.30 & 50 & 0.04 & 0.08 \\
\tabsp
(500, 0.25, 0.5) & 0  & 24.89 & 600.00 & 600.00 &  0 &  5.23 & 600.00 & 600.00 & 50 & 0.08 & 0.51 \\
(500, 0.25, 1.0) & 0  & 20.75 & 600.00 & 600.00 &  0 & 10.99 & 600.00 & 600.00 & 50 & 0.04 & 0.19 \\
(500, 0.75, 0.5) & 0  &  0.49 & 600.00 & 600.00 & 49 &  0.00 &  34.28 & 600.00 & 50 & 0.61 & 0.74 \\
(500, 0.75, 1.0) & 0  &  1.35 & 600.00 & 600.00 &  2 &  0.38 & 584.65 & 600.00 & 50 & 0.13 & 0.18 \\
\tabsp
(1000, 0.25, 0.5) & 0  & 28.40 & 600.00 & 600.00 &  0 &  6.27 & 600.00 & 600.00 & 50 & 0.33 & 1.96 \\
(1000, 0.25, 1.0) & 0  & 26.72 & 600.00 & 600.00 &  0 & 14.38 & 600.00 & 600.00 & 50 & 0.16 & 1.10 \\
(1000, 0.75, 0.5) & 0  &  0.79 & 600.00 & 600.00 &  1 &  0.03 & 598.81 & 600.00 & 50 & 1.83 & 2.27 \\
(1000, 0.75, 1.0) & 0  &  2.16 & 600.00 & 600.00 &  0 &  0.62 & 600.00 & 600.00 & 50 & 0.37 & 0.54 \\
\midrule
Sum & 302 & \multicolumn{3}{c}{} & 354 & \multicolumn{3}{c}{} & 800 \\	
Avg & & 7.59 & 391.60 & & & 2.78 & 344.85 & & & 0.24 & \\
Max & & 28.40 & & 600.00 & & 14.38 & & 600.00 & & & 2.27 \\
\bottomrule
\end{tabular}
\end{center}
\end{small}
\end{table}

\section{Assortment Optimization under Mixtures of Multinomial Logits} \label{sec:mmnl}

{As already mentioned, the \probname appears as a subproblem in a number of applications such as the assortment optimization problem under the MMNL choice model. In this section, we show (1) how our solution approach can be adapted to this setting and (2) that it provides a significant computational boost with respect to exact methods from the literature. The MMNL choice model assumes that a customer visiting a store belongs with a certain probability $\alpha^g$ to a certain customer type $g \in \mathcal{G}=\{1,....,G\}$ with preference weights $\pmb{v}^g=v_1^g,..., v_n^g$. In this context, an assortment that maximizes the expected revenue can be identified by solving the optimization problem}

\begin{equation}
\label{eq:AOPMMNL}
    z^* = \max_{\pmb{x} \in \{0,1\}^n} \Bigg\{ \sum_{g \in \mathcal{G}} \alpha^g \frac{\sum_{j \in P} r_j v_j^gx_j}{v_0^g + v^g(\pmb{x})} \Bigg\}.
\end{equation}




{\cite{Bront2009} and  \cite{Mendez-Diaz2014} propose MILP formulations for solving \eqref{eq:AOPMMNL}. \cite{Feldman2015} suggest to compute tight upper bounds to \eqref{eq:AOPMMNL} by means of Lagrangian relaxation. 
Specifically, rather than optimizing the assortment over the whole population of customers, their approach optimizes the assortment for each customer type separately, penalizing differences among customer-specific assortments by the use of Lagrangian multipliers  $\lambda_j^g$ for $j \in P$ and $g \in \mathcal{G}$. For a given set of multipliers $\pmb{\lambda}$, the optimal assortment for customer type $g\in \mathcal{G}$ is thus obtained by solving the problem}

\begin{equation}
\label{eq:AOPMMNL_lagrag}
    z^g(\lambda) = \max_{\pmb{x} \in \{0,1\}^n} \Bigg\{\frac{\sum_{j \in P} r_j v_j^gx_j}{v_0^g + v^g(\pmb{x})} -\sum_{j \in P} \lambda^g_jx_j\Bigg\}.
\end{equation}

{Here, $\lambda_j^g <0$ penalizes the introduction of item $j\in P$ in the assortment of customer type $g\in \mathcal{G}$, while  $\lambda_j^g >0$ encourages having $j$ in the assortment. Then, one may compute the \textit{best} upper bound by finding a set of multiplier $\pmb{\lambda}^*$ such that}

\begin{equation}
\label{eq:lagrdual}
    z(\pmb{\lambda^*}) = \min_{\pmb{\lambda}} \sum_{g \in \mathcal{G}}\alpha^g z^g(\pmb{\lambda}).
\end{equation}

 {It is easy to notice that \eqref{eq:AOPMMNL_lagrag} is an \probname in which Lagrangian multipliers $\lambda^g_j$ correspond to product costs $c_j^g$ for customer type $g\in \mathcal{G}$.
 This suggests that our exact solution approach for the \probname can be used as a subroutine for solving problem \eqref{eq:AOPMMNL}, provided that we are able to find an optimal set of multipliers $\pmb{\lambda}^*$. In their work, \cite{Feldman2015} use an iterative procedure based on subgradient search to find a set of Lagrangian multipliers $\widehat{\pmb{\lambda}}$ 
 and an approximate method to calculate an upper bound to $z^g(\pmb{\lambda}^k)$ in each iteration $k=1, \dots, K$. While their method does not guarantee to derive an optimal set of multipliers, it generally provides good performances. In this work, we use the same subgradient search, but use LLRS (instead of an approximate method) to compute $z^g(\pmb{\lambda}^k)$ for customer type $g \in \mathcal{G}$ at each iteration $k=1, \dots, K$. The algorithm terminates either when a time limit is reached or when the step size of the subgradient search becomes smaller than $\num{1e-8}$.  
 %
 
 \smallskip
 The following three adaptations are made to LLRS when solving \eqref{eq:lagrdual}:}
{
\begin{enumerate}
    \item \final\ is adapted as formulation~\eqref{prob:milp} allows for invalid solutions in which $x_j=1$ (i.e., $j\in P$ is in the assortment) and $u_j=0$ (i.e., $j$ is never chosen) even if $v_j>0$, which may arise in the case of negative product costs. Thus, we add constraints 
    \begin{equation}
    \label{eq:new_ctr}
    v_0 u_j \geq v_j u_0 + (x_j - 1)v_j \quad \forall j\in P  
    \end{equation}
    forcing inequality~\eqref{prob:milp:prob:ratios} to equality when $x_j=1$. Note that if $x_j=0$, then $u_j=0$ holds due to constraints~\eqref{prob:milp:force}. Thus, \eqref{eq:new_ctr} resolves to $0 \geq v_j (u_0 - 1)$, which is always satisfied.
    \item A (high-quality) primal solution \eqref{eq:AOPMMNL} is derived in each iteration $k=1, \dots K$ of the subgradient search using the (optimal) assortments $\mathcal{X}^{g,k}\subseteq P$ derived by LLRS for customer type $g\in \mathcal{G}$ and current multipliers $\pmb{\lambda}^k$. Since these assortments may be different from one another, we simply use their union $\mathcal{X}^k=\bigcup_{g\in\mathcal{G}}\mathcal{X}^{g,k}$ to obtain a heuristic solution to \eqref{eq:AOPMMNL}.
    \item Solving \final\ is skipped at iterations of the subgradient search in which the Lagrangian multipliers are likely not optimal. Although \final\ can typically be solved in fractions of a second (cf.\ \cref{tab:addInfo}), its solution is much slower than the time needed for the enhanced bounding procedure. Since, at a given iteration $k$, any feasible solution can be used to compute a subgradient of $\pmb{\lambda}^k$, we use heuristic solution $\mathcal{X}^k$ instead of solving \final\ in most iterations. Ideally, we resort to \final\ only at the last few iterations, when the set of multipliers is (close-to-)optimal and solving \eqref{eq:AOPMMNL_lagrag} to optimality may thus be required to close the final optimality gap.
    Practically, we execute \final\ only when the primal bound $\un{z}^k$ and the upper bound $\ov{z}^k$ returned by the bounding procedure are relatively close, i.e., $\ov{z}^k - \un{z}^k < \epsilon$. In our experiments we used $\epsilon=0.01$.
    
\end{enumerate}
}

\subsection{Data generation}\label{datageneration}
{To test LLRS with the adaptations described above, we generate a set of test instances as described in \cite{Feldman2015}\footnote{To this end, we used a code provided by the authors of \cite{Feldman2015}.}. Specifically, the number of products is set to $n=100$ and their revenues are sampled uniformly at random from the interval $[0,2\,000]$. Probabilities $\alpha^g$ for $g \in \mathcal{G}$ are proportional to customer-specific quantities $\beta^g$ that are sampled uniformly at random from the interval $[0,1]$ so that $\alpha^g = \frac{\beta^g}{\sum_{c \in \mathcal{G}} \beta^{c}}$. For each customer type, the products are randomly partitioned into the 40 staple and 60 specialty products. 
Specialty products are meant to simulate products for which the preference weight is either particularly high or particularly low. In particular, the preference weight for a generic product $j \in P$ is proportional to the quantity $k_j X^g_j$,  where $X_j^g$ is uniformly distributed in $[0.3,0.7]$ for staple products, and in $[0.1,0.3]\cup[0.7,0.9]$ for specialty products, respectively. Parameter $k_j$, used to further vary the magnitude of preference weights associated to each product, is sampled from the uniform distribution in $[1,\ov{K}]$. Finally, for each customer type $g \in \mathcal{G}$, the probability $P_0^g$ is sampled from the uniform distribution in $[0, \ov{P}_0^g]$, and preference weights are set so that a customer of type $g\in \mathcal{G}$ leaves without any purchase with probability $P_0^g$ when all products are offered. This is obtained by setting $v_j^g = \frac{k_j X_j^g (1-P_0^g)}{P_0^g \sum_{i \in P} k_j X_j^g}$ for each product $j\in P$ and customer class $g\in \mathcal{G}$. Instances are generated by varying the number of customer classes $G \in \{25,50,75\}$ and the parameters $\ov{K}$ and $\ov{P}_0^g$ in $\{5,10,20\}$ and $\{0.6,0.8,1\}$, respectively. One hundred instances are generated for each combination of parameters, resulting in 2\,700 instances in total.
}


\vspace{0.1cm}

\subsection{Computational results}\label{comput_results_mmnl}
{\cref{tab:MMNL} summarizes the computational performance of the modified version of LLRS embedded in the subgradient search and compares to the MILP formulation by \citet{Bront2009} solved with Cplex 20.1. For both methods, we report the number of instances that are solved to optimality (\textsf{opt}), the average optimality gap (\textsf{gap}) in percentage, the average computing time (\cpuavg) and the maximum one (\cpumax) on each class of instances. As for the \probname, all experiments are conduced on a single core of a machine with 500GB-RAM and an Intel(R) Xeon(R)Gold 6142 with 2.60GHz CPU. A time limit of ten minutes is imposed on each experiment. For this set of experiments, we used $\rho^f=\num{1e-2}$ and $\rho^\ell=\num{1e-4}$. Compared to the experiments in \cref{sec:results}, a larger value of $\rho^\ell$ was chosen to further speed-up the bounding procedure that, in this context, is executed thousands of times (once per customer type and iteration of the subgradient search). This is done also considering that \final\ needs to be executed only at the last iterations, thus allowing to not suffer significantly from the looser bounds obtained with this parameter setting.
The results given in \cref{tab:MMNL} show that LLRS is able to solve 2\,206 out of 2\,700 instances (roughly 82\%) optimally with an average computing time of about 10 seconds. The average optimality gap is less than 0.005\% and, while not reported in the table, the maximum optimality gap is 0.52\%. These numbers confirm the high quality of the solutions returned by the procedure even when optimality could not be proven. In comparison, the MILP formulation by \citet{Bront2009} is able to solve  only 238 of the instances (roughly 9\%) to optimality within the timelimit of ten minutes. Furthermore, the average optimality gaps from that method are approximately ten times larger than those obtained by LLRS. Compared to the MILP formulation, our approach seems to better scale with the number of customer classes, running on instances with 75 customer classes with an average computing time of 18 seconds. We also performed a separate set of experiments in which coarse grids are not used to speed up more refined ones (achieved by setting $\rho^f=\rho^\ell=\num{1e-4}$),  whose results allow to assess the impact of the \textit{enhanced} bounding procedure. While we do not report the detailed results of these experiments for the sake of brevity, we note that the average computing times increased to roughly 34 seconds in that case, i.e., an increase of around 240\% relative to the average computing times of 10.31 seconds reported in \cref{tab:MMNL}. Similarly to the results discussed for the \probname (cf.\ \cref{tab:vsFT}), these results clearly confirm the importance of the enhanced bounding procedure. 
}

\begin{table}[!htpb]
    \centering\sffamily
    {
    \caption{Summary of the computational results of \our\ on the AOP under MMNL and comparison with the MILP formulation by \cite{Bront2009}}
    \label{tab:MMNL}    
    \begin{tabular}{crrrrrrrrrr}																		
    \toprule																		
		&	\multicolumn{4}{c}{MILP	by	\cite{Bront2009}}	&	\multicolumn{4}{c}{\our}	\\										
	\cmidrule(lr){2-5}\cmidrule(lr){6-9}																	
	$G,\ov{K},\ov{P}_0$	&	opt	&	gap	&	\cpuavg	&	\cpumax	&	opt	&	gap	&	\cpuavg	&	\cpumax	\\
        \midrule																		
	(25, ~5, 0.6)	&	22	&	0.0104	&	517.17	&	600.14	&	87	&	0.0026	&	2.51	&	26.96	\\
	(25, ~5, 0.8)	&	24	&	0.0101	&	503.90	&	597.55	&	88	&	0.0009	&	2.30	&	27.64	\\
	(25, ~5, 1.0)	&	36	&	0.0086	&	453.43	&	597.63	&	85	&	0.0016	&	4.09	&	97.23	\\
	(25, 10, 0.6)	&	15	&	0.0147	&	540.38	&	600.09	&	79	&	0.0051	&	2.70	&	18.60	\\
	(25, 10, 0.8)	&	17	&	0.0145	&	544.43	&	600.05	&	83	&	0.0058	&	3.73	&	52.37	\\
	(25, 10, 1.0)	&	27	&	0.0156	&	498.27	&	600.08	&	90	&	0.0060	&	2.87	&	14.61	\\
	(25, 20, 0.6)	&	14	&	0.0150	&	550.14	&	600.08	&	80	&	0.0083	&	3.57	&	74.22	\\
	(25, 20, 0.8)	&	22	&	0.0141	&	522.42	&	600.06	&	86	&	0.0015	&	4.34	&	81.96	\\
	(25, 20, 1.0)	&	25	&	0.0122	&	507.12	&	600.04	&	88	&	0.0006	&	5.82	&	77.92	\\
 \tabsp
	(50, ~5, 0.6)	&	2	&	0.0224	&	592.74	&	596.89	&	87	&	0.0020	&	4.67	&	50.43	\\
	(50, ~5, 0.8)	&	4	&	0.0222	&	589.37	&	596.97	&	85	&	0.0029	&	4.44	&	28.39	\\
	(50, ~5, 1.0)	&	12	&	0.0187	&	565.61	&	596.84	&	84	&	0.0021	&	8.79	&	143.49	\\
	(50, 10, 0.6)	&	2	&	0.0230	&	592.72	&	596.85	&	81	&	0.0020	&	7.87	&	60.41	\\
	(50, 10, 0.8)	&	1	&	0.0221	&	595.42	&	596.33	&	87	&	0.0010	&	9.05	&	114.42	\\
	(50, 10, 1.0)	&	4	&	0.0219	&	581.31	&	596.25	&	79	&	0.0009	&	13.48	&	295.64	\\
	(50, 20, 0.6)	&	1	&	0.0259	&	594.23	&	596.22	&	80	&	0.0065	&	10.69	&	211.74	\\
	(50, 20, 0.8)	&	3	&	0.0260	&	589.41	&	596.36	&	79	&	0.0039	&	10.52	&	263.41	\\
	(50, 20, 1.0)	&	0	&	0.0231	&	596.03	&	596.41	&	81	&	0.0041	&	10.89	&	190.54	\\
 \tabsp
	(75, ~5, 0.6)	&	0	&	0.0251	&	595.95	&	596.21	&	85	&	0.0031	&	8.91	&	103.59	\\
	(75, ~5, 0.8)	&	1	&	0.0246	&	595.01	&	596.09	&	85	&	0.0016	&	11.10	&	156.65	\\
	(75, ~5, 1.0)	&	5	&	0.0216	&	587.26	&	597.84	&	79	&	0.0023	&	16.24	&	254.90	\\
	(75, 10, 0.6)	&	0	&	0.0276	&	595.79	&	596.22	&	80	&	0.0037	&	14.99	&	224.43	\\
	(75, 10, 0.8)	&	0	&	0.0289	&	595.76	&	596.05	&	74	&	0.0073	&	23.94	&	300.52	\\
	(75, 10, 1.0)	&	1	&	0.0267	&	592.49	&	596.00	&	77	&	0.0072	&	20.91	&	308.68	\\
	(75, 20, 0.6)	&	0	&	0.0298	&	595.69	&	596.04	&	75	&	0.0020	&	13.54	&	165.18	\\
	(75, 20, 0.8)	&	0	&	0.0309	&	595.79	&	596.12	&	69	&	0.0047	&	26.35	&	425.05	\\
	(75, 20, 1.0)	&	0	&	0.0324	&	595.66	&	596.06	&	73	&	0.0032	&	29.94	&	456.31	\\
    \midrule																		
	Sum	&	238	&		&		&		&	2206	&		&		&		\\
	Avg	&		&	0.0210	&	566.06	&		&		&	0.0034	&	10.31	&		\\
	Max	&		&		&		&	600.14	&		&		&		&	456.31	\\
        \bottomrule																		
    \end{tabular}															
    }
\end{table}

\section{Conclusions} \label{sec:conclusions}
In this work, we have proposed an exact method for the assortment optimization problem with product costs when
customers choose according to a multinomial logit model. This problem is (practically) relevant for several reasons. Indeed, product costs emerge in a variety of real-life scenarios, representing operational costs incurred by the firm for offering a certain product. Furthermore, despite its limitations, the multinomial logit model is 
one of the most studied discrete choice models, both in academia and industry, due to its interpretability and to the fact that it can be efficiently estimated. Moreover, as observed in a number of other works, the assortment optimization problem with product costs appears as a subproblem when optimizing decisions in the context of assortment optimization and network revenue management problems with side constraints. \cite{Kunnumkal2009} showed that solving this problem to optimality is NP-hard, motivating a number of approximate approaches to obtain \refone{lower} and \refone{upper} bounds in a tractable manner. We have tackled this problem by proposing a bounding procedure, which builds upon the grid-based approximate approach of \cite{Feldman2015}. In particular, we have exploited \refone{lower} and \refone{upper} bounds obtained from coarser grids which are computationally cheap to compute and alleviate the computational burden of denser grids. As a result, the iterative algorithm proposed in this article is able to compute tight \refone{lower} and \refone{upper} bounds in a fraction of a second, on instances with up to one thousand products. Furthermore, such instances are solved to optimality in roughly two tenths of a second, on average, by exploiting the obtained bounds. 
{We have also shown that our procedure can be readily adapted to handle a cardinality constraint on the size of assortments, with a limited impact on the computing times.
As a final contribution, we described how LLRS can be embedded in the solution approach by \cite{Feldman2015} to optimize assortments under the MMNL choice model. The computational results have shown that the resulting method can solve most of the considered benchmark instances with up to 75 products and 20 customer types optimally within relatively short time and derived extremely tight lower and upper bounds for those instances that could not be solved optimally. In particular, this method is much faster than an alternative exact approach by \citet{Bront2009} based on solving a MILP formulation. Comparing the obtained results to the bounding procedure by \citet{Feldman2015}, we observed that one appealing aspect of the algorithm proposed in this work is that primal solutions and proofs of optimality (in many cases) come almost for free since the increase in average computing time is very small. Future work may include the development of exact algorithms for assortment optimization problems with further side constraints or for problems in which customers choose according to other variants of the multinomial logit model. {As shown by our results on the MMNL, the} algorithm proposed in this article can be a key ingredient of such methods as the assortment optimization problem with product costs appears as a subproblem in several assortment optimization and network revenue management problems.
}

\section*{Acknowledgements}
We are indebted to four anonymous referees for their careful reading and meaningful remarks that we believe have helped us improve the paper and show the applicability of our approach to different contexts. We would also like to warmly thank Huseyin Topaloglu for sharing with us the code for the assortment optimization problem under the MMNL choice model.

\bibliographystyle{abbrvnat}
\bibliography{references} 

\begin{appendices}
    \section{Tuning of the hyper-parameters $\rho_f$ and $\rho_l$}\label{ap:sensitivity}
    \blue{In this section, we investigate the impact of the hyper-parameters \rhof and \rhol on the computational performance of \our. Furthermore, we propose rule-of-thumbs to guide practitioners in setting these two hyper-parameters. The same instances described in Section \ref{subsec:inst} were used for this set of experiments.}
    
    \blue{The first parameter, \rhof, should be set so that the number of intervals (i.e., Linear Knapsack relaxations) to solve at the first iteration is relatively small. To understand how small this number should be, we refer the reader to Table \ref{tab:rhof}, where we report the average and maximum computational time needed to solve MILP+, the bounding procedure, and the AOPC in general, as a function of \rhof when \rhol$=\num{1e-7}$. For each \rhof we also report the average number of intervals \#\intavg\ obtained at the first iteration. Note that \rhof does not impact MILP+ at all since the quality of the bound is determined by \rhol\ only. Also, the time needed to execute the bounding procedure decreases as \rhof increases. This is due to the fact that (i) we exploit coarse grids to prune more refined ones, and (ii) even if high values of \rhof may provide lose bounds, they are extremely fast to compute. In general, we suggest setting \rhof so that the number of linear knapsack relaxations to solve at the first iteration is about 100 or smaller.
     \begin{table}[!htpb]
    \centering\sffamily
    \begin{tabular}{rrrrrrrr}																
        \toprule																
        \rhof	&	\multicolumn{2}{c}{\final}	&	\multicolumn{3}{c}{Bounding Procedure}		&	\multicolumn{2}{l}{Total}	\\	
        \cmidrule(lr){2-3}\cmidrule(lr){4-7}\cmidrule{7-8}	
        	&	\cpuavg	&	\cpumax	&	\cpuavg	&	\cpumax	&	\#\intavg	&	\cpuavg	&	\cpumax	\\	
        \midrule																
        	$\num{1e-7}$	&	0.05	&	1.44	&	18.78	&	81.19	&	8,369,102.6	&	18.83	&	81.40	\\
        	$\num{1e-6}$	&	0.04	&	1.81	&	1.96	&	8.35	&	836,912.80	&	2.00	&	10.04	\\
        	$\num{1e-5}$	&	0.04	&	1.76	&	0.23	&	1.19	&	83,693.74	&	0.27	&	2.91	\\
        	$\num{1e-4}$	&	0.04	&	1.80	&	0.05	&	0.44	&	8,371.27	&	0.10	&	2.20	\\
        	$\num{1e-3}$	&	0.04	&	1.80	&	0.03	&	0.25	&	838.87	&	0.07	&	2.01	\\
        	$\num{1e-2}$	&	0.04	&	1.73	&	0.02	&	0.17	&	86.00	&	0.06	&	1.89	\\
        	$\num{1e-1}$	&	0.04	&	1.44	&	0.02	&	0.15	&	11.00	&	0.06	&	1.55	\\
        	$\num{1e-0}$	&	0.04	&	1.45	&	0.02	&	0.15	&	2.88	&	0.06	&	1.56	\\
        \bottomrule																
    \end{tabular}	
    \caption{Sensitivity analsysis w.r.t. \rhof, with \rhol$=\num{1e-7}$}
    \label{tab:rhof}
\end{table}

     In Table \ref{tab:rhol}, we conducted a similar study to investigate the computational performance of \our~as a function of \rhol, with \rhof~fixed to $\num{1e-2}$. In this case, the average time needed to solve \final~decreases as \rhol~decreases, as one may expect. In particular, the maximum time needed to solve \final~settles at 2-3 seconds for \rhol$\leq\num{1e-4}$. This value corresponds to a bound improvement $< 1\%$ with respect to the previous iteration. Such criterion may thus be used as a general rule of thumb for choosing \rhol. Furthermore, the results we report in Section \ref{comput_results_mmnl} confirm that even in this setting (\rhof$=\num{1e-2}$ and \rhol$=\num{1e-4}$) our bounding procedure allows to obtain significant speed-up (about 70\%) with respect to a vanilla procedure with \rhof=\rhol=$\num{1e-4}$. 
     
     Generally, the results described above provide computational evidence about the robustness of \our~with respect to the choice of its only two hyper-parameters.}
     
    \begin{table}[!htpb]
    \centering\sffamily
    \begin{tabular}{rrrrrrrrr}																		
        \toprule																		
        \rhol	&	\multicolumn{2}{c}{\final}	&	\multicolumn{2}{c}{Bounding Procedure}	&	\multicolumn{2}{c}{Total}	\\								\cmidrule(lr){2-3}\cmidrule(lr){4-5}\cmidrule{6-7}	
        	&	\cpuavg	&	\cpumax	&	\cpuavg	&	\cpumax	&		\cpuavg	&	\cpumax	\\	
        \midrule																		
        	$\num{1e-8}$	&	0.03	&	1.56	&	0.1164	&	1.6	&		0.15	&	3.16	\\
        	$\num{1e-7}$	&	0.04	&	1.73	&	0.0260	&	0.23	&		0.07	&	1.95	\\
        	$\num{1e-6}$	&	0.09	&	2.40	&	0.0081	&	0.07	&		0.10	&	2.42	\\
        	$\num{1e-5}$	&	0.31	&	3.99	&	0.0036	&	0.03	&		0.31	&	4.00	\\
        	$\num{1e-4}$	&	0.45	&	2.82	&	0.0018	&	0.03	&		0.45	&	2.82	\\
        	$\num{1e-3}$	&	0.73	&	39.50	&	0.0014	&	0.02	&		0.74	&	39.51	\\
        	$\num{1e-2}$	&	44.71	&	600.32	&	0.0010	&	0.02	&		44.71	&	600.33	\\
        \bottomrule																		
    \end{tabular}					
    \caption{Sensitivity analsysis w.r.t. \rhol, with \rhof$=\num{1e-2}$}
    \label{tab:rhol}
\end{table}
\end{appendices}
\end{document}